\documentclass[12pt]{article}
\usepackage{amssymb}
\usepackage{amsmath}
\usepackage{amsthm}
\usepackage{latexsym}
\usepackage{amsfonts}
\usepackage{graphicx}
\usepackage{graphics}

\theoremstyle{plain}
\newtheorem{thm}{Theorem}[section]   
\newtheorem{prop}{Proposition}[section]
\newtheorem{lem}{Lemma}[section]
\newtheorem{cor}{Corollary}[section]
\newtheorem{Def}{Definition}[section]

\theoremstyle{definition}

\newcommand{\dis}{\displaystyle}

\newcommand{\cp}{\mathcal{P}}

\newcommand{\ch}{\mathcal{H}}

\newcommand{\ld}{\ldots}

\newcommand{\el}{\ell}

\newcommand{\ra}{\rightarrow}

\newcommand{\bi}{\beta}
\newcommand{\ga}{\gamma }

\newcommand{\de}{\delta }

\newcommand{\vPsi} {{\varPsi}}
\newcommand{\e}{\varepsilon }
\newcommand{\f}{\varphi}
\newcommand{\Fi}{\varPhi}

\newcommand{\thi}{\theta }

\newcommand{\Thi} {{\mit\Theta}}

\newcommand{\la}{\lambda }
\newcommand{\mi}{\mu }

\newcommand{\si}{\sigma }
\newcommand{\ups}{\upsilon}
\newcommand{\ti}{\tau }

\newcommand{\oo}{\omega}

\newcommand{\R}{\mathbb{R}}

\newcommand{\Z}{\mathbb{Z}}
\newcommand{\N}{\mathbb{N}}

\newcommand{\bs}{\backslash}
\newcommand{\tsupp}{\mbox{supp}}
\newcommand{\ssum}{\sum\limits}
\newcommand{\bcup}{\bigcup\limits}

\newcommand{\oD}{\overline{D}}
\newcommand{\oB}{\overline{B}}
\begin{document}
\noindent
{\bf\Large{A new ergodic proof of a theorem of W. Veech}}
\footnotetext{\hspace{-0.5cm}P. Georgopoulos}\footnotetext{\hspace{-0.5cm}Department of Mathematics,
University of Athens, 15784, Athens,
Greece}\footnotetext{\hspace{-0.5cm}\vspace*{0.5cm}email:
pangeorgopoul@gmail.com \& pgeorgop@math.uoa.gr}
\vspace*{0.5cm}\\
%
{\bf Panagiotis Georgopoulos}
\vspace*{2cm}\\
%
%
\noindent
{\bf\footnotesize Abstract.} {\footnotesize Our goal in the present paper is to give a new ergodic proof of a well-known Veech's result, build upon our previous works [4,5].}\vspace*{0.2cm} \vspace*{0.2cm} 
\section{\normalsize Introduction}  
\noindent

W. Veech in his remarkable paper [11, Theorem 3] (see also [7, p. 235] and [8, Commentary of Problem 116, p. 203]), proved the following:

``Almost all'' sequences $(r_1,\ld,r_n,\ld)$ of positive integers have the following ``universal'' property: Whenever $G$ is a compact separable group and $z_1,z_2,\ld,z_n,\ld$ a sequence of elements of $G$ that generates a dense subgroup of $G$, then the sequence
$y_1,y_2,\ld,y_n,\ld$, where $y_n:=z_{r_1}\cdot z_{r_2}\ld z_{r_n}$ is uniformly distributed for the Haar measure on G. Veech called such sequences, ``uniformly distributed sequence generators''.

In \cite{5} we prove that:

``Almost all'' sequences
$(r_1,\ld,r_n,\ld)$ of positive integers have the following
``universal'' property: Whenever $(X,\mi)$ is a Borel probability measure,
compact metric space and $\Fi_1,\Fi_2,\ld,\Fi_n,\ld$ a sequence of
continuous, measure preserving maps on $(X,\mi)$, such that
the action (by composition) on $(X,\mi)$ of the semigroup with
generators $\Fi_1,\ld,\Fi_n,\ld$ is amenable (as discrete), uniquely ergodic and
non-sensitive on $\text{supp}\mi$, then for every $x\in X$ the sequence
$w_1,w_2,\ld,w_n,\ld$ where
\[
w_n:=\Fi_{r_n}(\Fi_{r_{n-1}}(\ld(\Fi_{r_2}(\Fi_{r_1}(x)))\ld))
\]
is uniformly distributed for $\mi$.

In the present paper we prove the next most special, albeit not direct, corollary of \cite{5}.

``Almost all'' sequences $(r_1,\ld,r_n,\ld)$ of positive integers have the following ``universal'' property: Whenever $G$ is a locally compact, amenable, separable group acting (continuously) on $(X,\mi)$ (a Borel probability measure compact metric space), by measure preserving homeomorphisms, such that the action is uniquely ergodic for $\mi$ and non-sensitive on $\tsupp\mi$ (it turns out that such an action is necessarily equicontinuous) and if $\Fi_n$, $n\in\N$ is a sequence in $G$ that generates (by composition) a dense semigroup in $G$ and $x\in X$, then the sequence $w_n:=\Fi_{r_n}(\Fi_{r_{n-1}}(\ld(\Fi_{r_2}(\Fi_{r_1}(x)))\ld))$, $n\in\N$ is uniformly distributed for $\mi$.\\
This completes investigation of [4,5] and gives Veech's theorem, at least for metrizable groups.

The new element in the present paper is Proposition \ref{prop4.1} that allows us to use a combination of the methods of [4,5]. In fact, in many aspects, most parts of the arguments of [4,5] are much simpler.

Next, let us explain how Veech's theorem falls in the frame of the above result.

Clearly, $G$ acts on $G$ (uniformly equicontinuously) by multiplication, i.e. for $g\in G$, $x\in G$, $(x,g)\mapsto x\cdot g$, $G$ is amenable (as compact) and the Haar measure $m_G$ is the unique invariant measure for this action. Also, the assumption that $z_1,z_2,\ld,z_n,\ld$ generate a dense subgroup of $G$, implies that the action of this subgroup on $G$ (by right translations) is uniquely ergodic for $m_G$.

On the other hand, the assumption that $z_1,z_2,\ld,z_n,\ld$ generate a dense subgroup of $G$, is equivalent to the assumption that $z_1,z_2,\ld,z_n,\ld$ generate a dense semigroup in $G$ (see [6, Theorem 9.16]).

Under these circumstances for $G$ metrizable, in view of our result (in particular for $x=e$) the sequence $y_n:=z_{r_1}\cdot z_{r_1}\ld z_{r_n}$, $n\in\N$ is uniformly distributed for $G$.

And a final remark: The general case, where the group $G$ is not necessarily metrizable, can be treated by similar methods, since the topology of $G$ is defined by a family of pseudometrics (see [3, Chapter IX, Section 11]).
\section{\normalsize The main results}  
\noindent

Throughout this paper $(p_1,\ld,p_n,\ld)$ is a probability sequence
with non-zero entries (i.e. $p_n>0$ for each $n$ and
$\ssum^\infty_{n=1}p_n=1$). We consider now the set of natural
numbers $\N=\{1,2,\ld\}$ endowed with the discrete topology. Then,
we take the one-point compactification of $\N$ and we get the
compact space $\widetilde{\N}:=\N\cup\{\infty\}$. Let
$(\widetilde{\N},m)$ be the measure space, where $m$ is a probability
measure on $\widetilde{\N}$, defined by $m(\{n\})\!\!=\!\!p_n$, for every
point $n$ on $\N$ and $m(\{\infty\})\!=\!0$. On the space
$Y:=\widetilde{\N}^{\Z}$, $\Z$ the integers, we consider the product measure
$\la:=\!\prod\limits^{+\infty}_{-\infty}\!m$ and the two-sided Bernoulli
shift $T\!:\!Y\!\ra\!Y$, with
$T(\{x_n\})\!=\!\{y_n\}$, where $y_n=x_{n+1}$, for every $n\in\Z$.

Also, throughout this paper, $G$ is an amenable, locally compact separable group acting (continuously) on a Borel probability measure, compact metric space $(X,\mi)$ and the action is uniquely ergodic for $\mi$ and non-sensitive on $\tsupp\mi$. It turns out (see Corollary \ref{cor4.1}), that such an action is necessarily equicontinuous.

Next, let $\Fi_1,\ld,\Fi_n,\ld$ be a sequence in $G$, that generates a dense semigroup in $G$. (Note that the action of this semigroup in $(X,\mi)$ is also uniquely ergodic).

We set up the skew product
\[
\vPsi:X\times Y\ra X\times Y \ \ \text{defined by} \ \
\vPsi(x,r):=(\Fi_{r_1}(x),T(r))
\]
where $r:=(\ld,r_{-n},\ld,r_{-1},r_0,r_1,\ld,r_n,\ld)$, conventionally we set\\
\[
\Fi_\infty\equiv Id_X \ \ (Id_X \; \text{the identity on}\; X).
\]

Clearly $\vPsi$ is Borel measurable and $\mi\times\la$ is invariant under $\vPsi$.
\begin{thm}\label{thm2.1}
If $\ti$ is a Borel probability measure on $X\times Y$, invariant
for $\vPsi$, such that the projection of $\ti$ on $Y$ equals
$\la$, then $\ti$ coincides with $\mi\times\la$.
\end{thm}

From the above theorem, taking $r=(\ld,r_{-n},\ld,r_{-1},r_0,r_1,\ld,\linebreak r_n,\ld)\in\N^\Z$ a generic point for $T$, it is easily seen, using some standard results (see [5, pp.\ 193-194]), that $(r_1,\ld,r_n,\ld)$ has the property mentioned in the abstract.
\section{\normalsize Invariant measures for continuous maps}  
\noindent

The space $M(X)$ of all Borel probability measures on $X$ is
metrizable in the weak$^\ast$ topology. If
$\big\{f_n\big\}^\infty_{n=1}$ is a dense subset of $C(X)$ (the
space of continuous functions on $X$), then
\[
d(\si,\nu):=\sum^\infty_{n=1}\frac{|\int f_nd\si-\int
f_nd\nu|}{2^n\|f_n\|}
\]
is a metric on $M(X)$ giving the weak$^\ast$ topology. Also, $M(X)$
is compact in this topology.

For $\Fi:X\ra X$ continuous, hence Borel measurable,
we have the continuous affine map
\[
\f:M(X)\ra M(X) \ \ \text{given by} \ \
(\f\si)(B)=\si(\Fi^{-1}(B))
\]
for $B$ a Borel set.

We have
\begin{thm}\label{thm3.1}
Let $F_m$, $m\in\N$ be a F\"{o}lner sequence in $G$. For $\nu\in M(X)$ and $m\in\N$ we consider the measures
\[
\mi^\nu_m:=\frac{1}{m_G(F_m)}\int_{F_m}\f(\nu)\,dm_G(\Fi)
\]
(where $m_G$ is the Haar measure on $G$), or more concretely
\[
\int_Xf(x)\,d\mi^\nu_m(x):=\frac{1}{m_G(F_m)}\int_{F_m}\int_Xf(\Fi(x))\,d\nu(x)\, dm_G(\Fi)
\]
for every $f\in C(X)$ and every $m\in\N$.\\
Then, $d(\mi^\nu_m,\mi)\ra0$ for $m\ra\infty$ uniformly for $\nu\in M(X)$.
\end{thm}
\noindent
%
{\em Proof}. Suppose that the conclusion of the theorem does not hold. Then, there exist an $\e>0$, a subsequence $F_{m_n}$, $n\in\N$ of $F_m$, $m\in\N$ and a sequence $\nu_n$, $n\in\N$ in $M(X)$ such that
\begin{eqnarray}
d(\mi^{\nu_n}_{m_n},\mi)>\e.   \label{eq1}
\end{eqnarray}
For $f\in C(X)$ we have
\[
\int_Xf(x)\,d\mi^{\nu_n}_{m_n}(x):=\frac{1}{m_G(F_{m_n})}\int_{F_{m_n}}\int_Xf
(\Fi(x))\,d\nu_n(x)\,dm_G(\Fi)
\]
and for $H\in G$ $(h:M(X)\ra M(X)$ the induced map),
\begin{align*}
\int_Xf(x)\,dh(\mi^{\nu_n}_{m_n}(x)):&=\frac{1}{m_G(F_{m_n})}\int_{F_{m_n}}\int_X
f(H\circ\Fi(x))\,d\nu_n(x)\,dm_G(\Fi)\\
&=\frac{1}{m_G(F_{m_n})}\int_{H\,F_{m_n}}\int_Xf(\Fi(x))\,d\nu_n(x)\,dm_G(\Fi).
\end{align*}
So
\begin{align*}
\bigg|\int_Xf(x)\,d\mi^{\nu_n}_{m_n}(x)&-\int_Xf(x)\,dh(\mi^{\nu_n}_{m_n}(x))\bigg|\\
&\le
\frac{1}{m_G(F_{m_n})}\int_{F_{m_n}\triangle HF_{m_n}}\int_X|f(\Fi(x)|d\nu_n(x)\,dm_G(\Fi)\\
&\le\frac{m_G(F_{m_n}\triangle HF_{m_n})}{m_G(F_{m_n})}\|f\|_\infty\ra0 \ \ \text{for} \ \ n\ra\infty.
\end{align*}

Hence, every $w^\ast$-limit of the sequence $\mi^{\nu_n}_{m_n}$, $n\in\N$ is invariant under the action of $G$, so equals $\mi$ contradicting (\ref{eq1}). \hfill$\square$
\section{\normalsize Some results on amenable, non-sensitive actions}\label{sec4}
\noindent

We recall the following
\begin{Def}\label{Def4.1}
(See also [1, p. 23]) A continuous action of a group $G$, on a compact metric space $(X,\ups)$ ($\ups$ denotes the metric on $X$), is called sensitive on a subset $X'\subset X$, if there exists a $\bi>0$, such that for every $x\in X'$ and $\de>0$, there exist a $y\in X$ with $\ups(x,y)<\de$ and an $h\in G$, such that $\ups(h(x),h(y))\ge\bi$. Otherwise the action is called non-sensitive on $X'\subset X$.
\end{Def}

We set for $k\in\N$

$E_k:=\{x\in X:$ there exists an open neighborhood $U$ of $x$ such that  \\
\hspace*{3.6cm}$x_1,x_2\in U\Rightarrow\ups(\Fi(x_1),\Fi(x_2))<\dfrac{1}{k}$, for all $\Fi\in G\}$.\smallskip\\
Clearly, $E_k$ is open and since the action of $G$ is non-sensitive on $\text{supp}\mi$, $E_k\cap\text{supp}\mi\neq\emptyset$, for every $k\in\N$.

Note that a $x\in X$ is an equicontinuity point for $G$, if for every $\e>0$ there exists a $\de>0$ such that $\ups(x,y)<\de$ implies $\ups(\Fi(x),\Fi(y))<\e$, for every $\Fi\in G$. Clearly, $\dis\bigcap^\infty_{k=1}E_k$ is the set of equicontinuity points for $G$.
\begin{lem}\label{lem4.1}
Let $k\in\N$. Then for every $x\in X\backslash E_k$ there exists a $\Fi_{i_x}\in G$, such that $\Fi_{i_x}(x)\in E_k$.
\end{lem}
\noindent
{\em Proof}. For $k\in\N$, the set
\[
Q_k:=(X\backslash E_k)\Big\backslash \bigcup_{\Fi\in G}\Fi^{-1}(E_k)
\]
is compact and forward invariant under the elements of $G$.

In case that $Q_k\neq\emptyset$, by an application of Day's fixed point theorem [2, Theorem 1], there exists a Borel probability measure $\ti$ supported on $Q_k$ and invariant under $G$, so $\ti=\mi$. But this contradicts the fact that $E_k\cap\text{supp}\mi\neq\emptyset$, for every $k$. So, $Q_k=\emptyset$
and the conclusion of the lemma follows\linebreak immediately. \hfill$\square$
\begin{cor}\label{cor4.1}
The group $G$ acts on $X$ equicontinuously.
\end{cor}
{\em Proof}. Since the maps $\Fi:X\ra X$, $\Fi\in G$ are open (as homeomorphisms), it is easily seen that $\Fi(E_k)\subseteq E_k$ for every $k\in\N$ and $\Fi\in G$.

Let $x\in X$. Suppose, if possible, that $x$ is not an equicontinuity point for the action of $G$ in $X$. Then
\[
x\in X\big\backslash \bigcap^\infty_{k=1}E_k.
\]
So, there exists a $k_0\in\N$ such that $x\notin E_{k_0}$. By the previous lemma, there exists a $\Fi_{i_x}\in G$ such that $\Fi_{i_x}(x)\in E_{k_0}$.
Since $\Fi(E_{k_0})\subseteq E_{k_0}$, for every $\Fi\in G$, clearly we have
$\Fi^{-1}_{i_x}\circ\Fi_{i_x}(x)=x\in E_{k_0}$, a contradiction. \hfill$\square$\vspace*{0.2cm}

We set $Seq:=\bigcup\limits^\infty_{n=1}\N^n$ the set of finite sequences of positive integers, and for $r=(r_1,\ld,r_n)\in Seq$, $\Fi_r:=\Fi_{r_n}\circ\cdots\circ\Fi_{r_1}$, $\f_r:=\f_{r_n}\circ\cdots\circ\f_{r_1}$ and $\Thi:=\{\f_r:r\in Seq\}$.

Under the above setting we have the following proposition, which is the new element that gives the possibility to use a combination of the methods of [4,5] in the present situation (see [5, Proposition 3.1]).
\begin{prop}\label{prop4.1}
There exists a sequence $\rho_m$, $m\in\N$ in $conv(\Thi)$ (the convex hull of $\Thi$) such that
\[
d(\rho_m(\si),\mi)\ra0 \quad \text{uniformly for} \quad \si\in M(X).
\]
\end{prop}
{\em Proof}. By Theorem \ref{thm3.1}, we can assume that there exist a F\"{o}lner sequence $F_m$, $m\in\N$ in $G$, and $\e_m>0$, $m\in\N$ with $\e_m\ra0$ for $m\ra\infty$ such that setting, for $\si\in M(X)$, $\mi^\si_m\in M(X)$ with
\[
\int_Xfd\mi^\si_m:=\frac{1}{m_G(F_m)}\int_{F_m}\int_Xf(\Fi(x))d\si dm_G(\Fi) \quad \text{for} \quad f\in C(X)
\]
we have
\begin{eqnarray}
d(\mi^\si_m,\mi)<\e_m \quad\text{for}\quad m=1,2,\ld \quad\text{and} \quad \si\in M(X). \label{eq2}
\end{eqnarray}

Let $D\!\subseteq\! X$ be denumerable, with $\oD\!=\!X$. We enumerate $D=\{x_i:i\in\N\}$ and set $A:=\{\de_{x_i}:x_i\in D$, $i\in\N$ and $\de_{x_i}$ is the Dirac measure on $x_i\}$ $(\subseteq M(X))$.

Also, let $\{f_n:n\in\N\}(\subseteq C(X))$ be dense in $C(X)$ (clearly $\{f_n:n\in\N\}$ defines the metric on $M(X)$, see above).

Let $m\in\N$. For $n=1,\ld,m$, $i=1,\ld,m$ we set
\[
g^i_n:G\ra\R, \quad \text{where} \quad g^i_n(\Fi)=\int_Xf_n\circ\Fi(x)d\de_{x_i}.
\]
It is easily seen, that the above $g^i_n$ are continuous.

Clearly, for $m\in\N$ and $n=1,\ld,m$, $i=1,\ld,m$ we have
\begin{eqnarray}
\int_X f_nd\mi^{\de_{x_i}}_m=\frac{1}{m_G(F_m)}\int_{F_m}g^i_n(\Fi)dm_G.  \label{eq3}
\end{eqnarray}
We set $B:=\{\Fi_\el:\el\in Seq\}$. By assumption we have $\oB=G$.\\
By [9, Chapter II, Theorem 6.3], for $m\in\N$ there exists a convex combination
\[
\sum^{k_m}_{k=1}\la_k\de_{\Fi_{\el_k}}, \quad \Fi_{\el_k}\in B, \quad k=1,\ld,k_m
\]
of Dirac measures on $M(G)$, such that for $i=1,\ld,m$ and $n=1,\ld,m$
\[
\bigg|\frac{1}{m_G(F_m)}\int_{F_m}g^i_n(\Fi)dm_G-\sum^{k_m}_{k=1}\la_kg^i_n(\Fi_{\el_k})\bigg|\le\e_m
\cdot\|f_n\|.
\]
So, in view of (\ref{eq3}) and the definition of the $g^i_n$'s, for $m\in\N$, $i=1,\ld,m$ and $n=1,\ld,m$
\begin{eqnarray}
\bigg|\int_Xf_nd\mi^{\de_{x_i}}_m-\sum^{k_m}_{k=1}\la_k
\int_Xf_n\circ\Fi_{\el_k}(y)d\de_{x_i}\bigg|\le\e_m\cdot\|f_n\|.  \label{eq4}
\end{eqnarray}

Setting $\rho_m:=\ssum^{k_m}_{k=1}\la_k\f_{\el_k}$, we have for $m\in\N$, $i=1,\ld,m$ and $n=1,\ld,m$
\[
\bigg|\int_Xf_nd\mi^{\de_{x_i}}_m-\int_Xf_nd\rho_m(\de_{x_i})\bigg|\le\e_m
\cdot\|f_n\|.
\]
So, for $m\in\N$ and $i=1,\ld,m$
\begin{align}
d(\mi^{\de_{x_i}}_m,\rho_m(\de_{x_i}))&\le \e_m\bigg(1-\frac{1}{2^m}\bigg)+2\sum^\infty_{n=m+1}\frac{1}{2^n}
 \nonumber\\
 &<\e_m+\frac{1}{2^{m-1}}.  \label{eq5}
\end{align}

Combining (\ref{eq2}) and (\ref{eq5}), it follows that for $m\in\N$ and $i=1,\ld,m$
\begin{eqnarray}
d(\rho_m(\de_{x_i}),\mi)<2\e_m+\frac{1}{2^{m-1}}.  \label{eq6}
\end{eqnarray}
\noindent
{\bf Claim 1.} $\rho_m(\de_x)\ra\mi$ uniformly for $x\in X$.

Let $\e>0$. There exists an $m_0\in\N$ such that
\[
\frac{1}{2^{m-1}}<\e \quad \text{and} \quad \e_m<\e \quad \text{for} \quad m>m_0.
\]
Let $f_1,\ld,f_{m_0}$. For the given $\e>0$ there exists $\de>0$, such that for $x,x'\in X$ with $v(x,x')<\de$
\[
|f_n(x)-f_n(x')|<\e\cdot\|f_n\| \quad \text{for} \quad n=1,\ld,m_0
\]
(where $v$ denotes the metric on $X$).\\
Since $B:=\{\Fi_\el:\el\in Seq\}$ is equicontinuous, for the above $\de>0$ there exists $\thi>0$ such that for $y,y'\in X$ with $v(y,y')<\thi$
\[
v(\Fi_\el(y),\Fi_\el(y'))<\de \quad \text{for every} \ \ \Fi_\el\in B.
\]

Since $\oD=X$, there exists an $m_\ast>m_0$ such that for every $x\in X$, there exists a $x_{i_\ast}\in D$, $i_\ast\in\{1,\ld,m_\ast\}$ with $v(x_{i_\ast},x)<\thi$.

So, for every $x\in X$, $m>m_\ast$ and $n=1,2,\ld,m_0$ we have
\[
\bigg|\sum^{k_m}_{k=1}\la_k\int_Xf_n\circ\Fi_{\el_k}(y)d\de_{x_{i_\ast}}-
\sum^{k_m}_{k=1}\la_k\int_Xf_n\circ\Fi_{\el_k}(y)d\de_x\bigg|<\e\cdot\|f_n\|
\]
and in view of (\ref{eq4}), since $i_\ast\in\{1,\ld,m_\ast\}$, we have for every $x\in X$, $m>m_\ast$ and $n=1,2,\ld,m_0$
\[
\bigg|\int_Xf_nd\mi^{\de_{x_{i_\ast}}}_m-\int_Xf_nd\rho_m(\de_x)\bigg|\le2\cdot\e
\cdot\|f_n\|
\]
(note that $\e_m<\e$ for $m>m_\ast>m_0$).

So, for every $x\in X$, $m>m_\ast$ we have
\[
d(\mi^{\de_{x_{i_\ast}}}_m,\rho_m(\de_x))<2\e+\frac{1}{2^{m_0-1}}.
\]
Finally, by (\ref{eq6}) we have that for every $x\in X$ and $m>m_\ast$
\[
d(\rho_m(\de_x),\mi)<\bigg(2\e_m+\frac{1}{2^{m-1}}\bigg)+\bigg(2\e+\frac{1}{2^{m_0-1}}\bigg)
<4\e+2\e=6\e
\]
(note that for $m>m_\ast>m_0$, $\e_m<\e$ and $\dfrac{1}{2^{m-1}}<\e$).\vspace*{0.2cm}\\
\noindent
{\bf Claim 2.} $\rho_m(\si)\ra\mi$ uniformly for $\si\in\Big\{\ssum^s_{k=1}\la_k\de_{x_k}:\ssum^s_{k=1}\la_k=1$, $x_k\in D\Big\}$.

Indeed, the claim holds from Claim 1, since $\rho_m(\si)$ is a convex combination of measures of the form $\rho_m(\de_x)$, $x\in X$.

Finally, $\rho_m(\si)\!\ra\!\mi$ uniformly for every $\si\in M(X)$, since the set $\Big\{\ssum^s_{k=1}\la_k\de_{x_k}\!:\ssum^s_{k=1}\la_k=1$, $x_k\in D\Big\}$ is dense in $M(X)$ by [9, Chapter II, Theorem 6.3].\hfill$\square$\vspace*{0.2cm}

The following lemma is a simplification of [5, Lemma 4.4].
\begin{lem}\label{lem4.2}
Let $\rho_m$, $m\in\N$ a sequence in $\text{conv}(\Thi)$ as in Proposition \ref{prop4.1}, $\nu_m$, $m\in\N$, $h_m$, $m\in\N$ sequences in $M(X)$ and $Seq$ respectively and $f\in C(X)$. Then
\[
\int_X f\circ\Fi_{h_{m_\el}}d\rho_{m_\el}(\nu_{m_\el})\longrightarrow
\int_X fd\mi \quad \text{for} \quad \el\ra\infty,
\]
for some subsequence $m_\el$, $\el\in\N$, of $m\in\N$.
\end{lem}
\noindent
{\em Proof}. Since the action of $G$ on $X$ is equicontinuous, the sequence $\Fi_{h_m}$, $m\in\N$ is equicontinuous for every sequence $h_m$, $m\in\N$ in Seq. Then $f\circ\Fi_{h_m}$, $m\in\N$
 is equicontinuous, so by Arzela-Ascoli theorem it has a uniformly convergent subsequence
\[
f\circ\Fi_{h_{m_\el}}\overset{u}{\longrightarrow}\widetilde{f}\in C(X).
\]
Then for $\e>0$ there exists an $\el_1\in\N$ such that
\[
\|f\circ\Fi_{h_{m_\el}}-\widetilde{f}\|_\infty<\e \quad \text{for} \quad \el\ge\el_1.
\]
So
\begin{eqnarray}
\bigg|\int_X f\circ\Fi_{h_{m_\el}}d\rho_{m_\el}(\nu_{m_\el})-\int_X\widetilde{f}
d\rho_{m_\el}(\nu_{m_\el})\bigg|<\e \quad\text{for} \quad \el\ge\el_1.  \label{eq7}
\end{eqnarray}
On the other hand, by Proposition \ref{prop4.1} there exists an $\el_2\in\N$ such that
\begin{eqnarray}
\bigg|\int_X\widetilde{f}d\rho_{m_\el}(\nu_{m_\el})-\int_X\widetilde{f}d\mi\bigg|<\e \quad \text{for} \quad \el\ge\el_2.  \label{eq8}
\end{eqnarray}
By (\ref{eq7}) and (\ref{eq8}) there exists an $\el_0\in\N$ so that
\[
\bigg|\int_X f\circ\Fi_{h_{m_\el}}d\rho_{m_\el}(\nu_{m_\el})-\int_X\widetilde{f}d\mi
\bigg|<2\e \quad \text{for} \quad \el>\el_0.
\]
Hence
\[
\int_X f\circ\Fi_{h_{m_\el}}d\rho_{m_\el}(\nu_{m_\el})\longrightarrow\int_X\widetilde{f}
d\mi \quad \text{for} \quad \el\ra\infty.
\]
Now it suffices to show that $\dis\int_X fd\mi=\dis\int_X\widetilde{f}d\mi$.

Indeed, $\dis\int_X f\circ\Fi_{h_{m_\el}}d\mi=\dis\int_X fd\mi$, since the $\Fi_r$'s, $r\in Seq$ preserve the measure $\mi$ and $f\circ\Fi_{h_{m_\el}}\overset{u}{\longrightarrow}\widetilde{f}$, so $\dis\int_X fd\mi=\dis\int_X\widetilde{f}d\mi$. \hfill$\square$
\begin{cor}\label{cor4.2}
Let $\rho_m$, $m\in\N$, $\nu_m$, $m\in\N$, $h_m$, $m\in\N$ sequences as in Lemma \ref{lem4.2} and $K\subset X$ Jordan measurable, i.e. $\mi(\partial K)=0$ ($\partial K$ the boundary of $K$) with $\mi(K)>a$, for some $0<a<1$.
Then there exists an $m_{\el_0}\in\N$ such that
\[
\int_X\chi_K\circ\Fi_{h_{m_{\el_0}}}d\rho_{m_{\el_0}}
(\nu_{m_{\el_0}})>a.
\]
\end{cor}

The proof of the corollary is similar to that of [5, Corollary 4.3], so we omit it.
\section{\normalsize Some technical lemmata}\label{sec5}
\noindent

In the sequel, we assume the curriculum of notations and definitions of [4, Section \ref{sec5}]. For $A\subseteq\Z$, $pr_A:\N^\Z\ra\N^A$ denotes the natural projection and for $k\in\N$, $Z_k:=\{-k,\ld,0,\ld,k\}$.

We recall from \cite{4} and \cite{5} the following lemmata.
\begin{lem}\label{lem5.1}
Let $B\subseteq\N^{\Z}$ compact with $\la(B)>0$ and $\bi$ with $0<\bi<1$. Then there exists an $a=(a_{-k},\ld,a_{-1},a_0,a_1,\ld,a_k)\in\N^{\Z_k}$, for $k\in\N$ such that
\[
\frac{\la(pr^{-1}_{\Z_k}\{a\}\cap B)}{\la\big(pr^{-1}_{\Z_k}\{a\}\big)}>1-\bi.
\]
\end{lem}
\noindent
{\em Proof}. See [4, Lemma 5.1]. \hfill$\square$
\begin{lem}\label{lem5.2}
Let $F\subseteq Seq$ finite. Then there exists a $\bi$, $0<\bi<1$, such that, if $B\subseteq\N^\Z$ measurable, with $\la(B)>0$ and $a\in\N^{\Z_k}$ for some $k\in\N$ satisfying
\[
\frac{\la(pr^{-1}_{\Z_k}\{a\}\cap B)}{\la(pr^{-1}_{\Z_k}\{a\})}>1-\bi,
\]
then for sufficiently large $n$ $(n\ge n_1)$, there exists a $t_n\in\N^{n-2k-1}$ such that
\[
\la([\widetilde{pr}^{-1}\{(a,t_n,z,a)\}\cap T^{n+|z|}(B)]\cap[\widetilde{pr}^{-1}\{a\}\cap B])>0
\]
for all $z\in F$, (where $|z|$ denotes the length of $z$).
\end{lem}
\noindent
{\em Proof}. See [5, Lemma 6.1]. \hfill$\square$\medskip

The following lemma is highly technical and its meaning will be clear
in the proof of Theorem 6.2.
\begin{lem}\label{lem5.3}
Let $\nu$ be a Borel probability measure on $X\times Y$ singular with
respect to $\mi\times\la$, such that the projection of $\nu$ on
$Y$ coincides with $\la$. Then given $0<\oo<1$, $0<\thi<1$ and
$h:\R^+\ra\R^+$ a non-decreasing function, there exist $Q_k$,
$k=1,2,\ld,s$, $s\in\N$, disjoint compact subsets of $X$,
$K\subseteq X\bs\bigcup\limits^s_{k=1}Q_k$ compact, and
$B\subseteq Y$ compact, with $\la(B)>0$, such that
\end{lem}
\begin{enumerate}
\item[{\em (i)}] \; $\mi(K)>1-\oo$, $\mi(\partial K)=0$ $(\partial K$
{\em the boundary)}
\item[{\em (ii)}] \; setting $e:=$\,{\em distance
$\Big(K,\bigcup\limits^s_{k=1}Q_k\Big)>0$, we have}
\[
\text{{\em diameter}}\;(Q_k)<h(e) \ \ \text{for} \ \ k=1,2,\ld,s
\]
\item[{\em (iii)}] \;
$\nu_y\Big(\bigcup\limits^s_{k=1}Q_k\Big)>1-\thi$,\; for \; $y\in
B$
\item[{\em (iv)}] \; $|\nu_y(Q_k)-\nu_{y'}(Q_k)|<\dfrac{\thi}{s}$\; {\em for
every} \; $y,y'\in B$,\; $k=1,2,\ld,s$

{\em (where $\nu_y$ denotes the conditional measure induced by $\nu$ on
the fiber $X\times\{y\}$)}.
\end{enumerate}
\noindent
{\em Proof}. See [4, Lemma 6.1]. \hfill$\square$\vspace*{0.2cm} \\
\noindent
{\em Note}. Although the $\Fi$'s in \cite{4} are commutative, this is not used in the proof of [4, Lemma 6.1].\vspace*{0.2cm}

Under the assumptions of Lemma \ref{lem5.3}, we have the following
\begin{cor}\label{cor5.1}
Let $y_0\in B$, $B'\subset B$ measurable, with $\la(B')>0$ and $\cp\subset\{1,2,\ld,s\}$, such that
\[
\sum_{k\in\cp}\nu_{y_0}(Q_k)>1-\e, \ \ \text{for} \ \ 0<\e<1.
\]
Then
\[
\nu\bigg(\bigg(\bigcup_{k\in\cp}\overline{Q}_k\bigg)\times B'\bigg)>((1-\e)-\thi)\cdot\la(B').
\]
\end{cor}
\noindent
{\em Proof}. See [5, Corollary 5.1].\hfill$\square$
\section{\normalsize The proof of Theorem 2.1}\label{sec6}
\noindent

The proof of Theorem \ref{thm2.1} will be given in two major steps. First, we shall prove that if $\ti$ is absolutely continuous with respect to $\mi\times\la$ then $\ti$ coincides with $\mi\times\la$. Second, we shall prove that $\ti$ has a trivial singular part with respect to $\mi\times\la$. These two steps are described in Theorems \ref{thm6.1} and \ref{thm6.2}, respectively.

We have
\begin{thm}\label{thm6.1}
The measure $\mi\times\la$ is the unique Borel probability measure on $X\times Y$, invariant under $\vPsi$ and absolutely continuous with respect to $\mi\times\la$.
\end{thm}
\noindent
{\em Proof}. This follows from the ergodicity of the skew product $\vPsi$, see the random ergodic theorem in \cite{10}. \hfill$\square$\vspace*{0.2cm}\\
\noindent
{\em Remark}. Note that the use of the random ergodic theorem of Ryll-Nardzewski (see \cite{10}) gives immediately Theorem \ref{thm6.1}, so we can omit the lengthy proof of the ``first step'' that appears in [4, Proposition 5.1] and [5, Theorem 6.1].\vspace*{0.2cm}

The proof of the following theorem is an amalgamation of the proofs of [4, Theorem 7.1] and [5, Theorem 7.1].
\begin{thm}\label{thm6.2}
Let $\nu$ be a Borel probability measure on $X\times Y$ singular
with respect to $\mi\times\la$, such that the projection of $\nu$
on $Y$ coincides with $\la$. Then $\nu$ is not invariant under
$\vPsi$.
\end{thm}
\noindent
{\em Proof}. Suppose that the conclusion of the theorem does not hold i.e. $\nu$ is invariant
for $\vPsi$.

Since the semigroup $\ch$ generated by $\Fi_1,\ld,\Fi_n,\ld$ acts
equicontinuously on $X$ (by Corollary \ref{cor4.1}), if $\rho$ denotes the metric on $X$, then
clearly there exists a non-decreasing $h:\R^+\ra\R^+$, such that
for every $f\in\ch$ and $x,y\in X$ with $\rho(x,y)<h(\de)$
$(\de>0)$, then $\rho(f(x),f(y))<\de$. Now given
$0<\oo<\dfrac{1}{100}$, $0<\thi<\dfrac{1}{100}$ and $h$ as above,
by Lemma \ref{lem5.3} there exist $Q_k$, $k=1,\ld,s$, disjoint
compact subsets of $X$, $K\subseteq X\bs\bcup^s_{k=1}Q_k$ compact
and $B_1\subseteq Y:=\widetilde{\N}^{\Z}$ compact with
$\la(B_1)>0$ satisfying conditions (i), (ii), (iii), (iv) of the
lemma, (with $B_1$ in place of $B$).

Let $B'_1:=B_1\cap\N^{\Z}$. Then $\la(B'_1)=\la(B_1)>0$ and
by the regularity of $\la$, there exists some compact $B\subseteq
B'_1$, such that $\la(B)>0$. The set $B$ satisfies the conditions
of Lemma \ref{lem5.3}

We consider $\rho_m$, $m\in\N$ a sequence in $conv(\Thi)$ as in Proposition \ref{prop4.1}. Since $\rho_m\in conv(\Thi)$, there exist a finite $F_m\subset Seq$ and $\thi_z(m)>0$ for $z\in F_m$, such that $\sum\limits_{z\in F_m}\thi_z(m)=1$ and $\rho_m=\sum\limits_{z\in F_m}\thi_z(m)\f_z$.

By Lemma \ref{lem5.2} for each $F_m$, $(m\in\N)$ there exists a $\bi_m$, $0<\bi_m<1$, satisfying the conclusion of that lemma.

Applying Lemma \ref{lem5.1} repeatedly, we find for each couple
\[
B,\bi_m \qquad m=1,2,\ld
\]
a $k_m\in\N$ and an $a^{(m)}=(a^{(m)}_{-k_m},\ld,a^{(m)}_0,\ld,a^{(m)}_{k_m})\in\N^{\Z_{k_m}}$ satisfying
\begin{eqnarray}
\frac{\la(B\cap pr^{-1}_{\Z_{k_m}}\{a^{(m)}\})}
{\la(pr^{-1}_{\Z_{k_m}}\{a^{(m)}\})}>1-\bi_m  \label{eq9}
\end{eqnarray}
for $m=1,2,\ld\;.$

Next, applying Lemma \ref{lem5.2} repeatedly, taking in view of (\ref{eq9}), we find for each quadruple
\[
F_m,\bi_m,B,\;a^{(m)}\in\N^{\Z_{k_m}} \quad \text{for some} \quad k_m\in\N, \quad m=1,2,\ld,
\]
an $n_m\in\N$ and a $t_{n_m}\in\N^{n_m-2k_m-1}$ such that, setting $t_{n_m}=t_m$ for brevity in the notation,
\begin{eqnarray}
\la([\widetilde{pr}^{-1}\{a^{(m)},t_m,z,a^{(m)}\}\cap T^{n_m+|z|}(B)]\cap[\widetilde{pr}^{-1}\{a^{(m)}\}\cap B])>0,  \label{eq10}
\end{eqnarray}
for all $z\in F_m$.

In the sequel we fix some $y_0\in B$ and set
\[
\ga_k:=\frac{\nu_{y_0}(Q_k)}{\nu_{y_0}\Big(\dis\bigcup^s_{i=1}Q_i\Big)}, \quad k=1,2,\ld,s.
\]
We fix $x_k\in Q_k$, $k=1,2,\ld,s$ and consider the probability measure
\[
\ti:=\sum^s_{k=1}\ga_k\de_{x_k}, \quad (\de_{x_k}\ \ \text{the Dirac measure)}.
\]

At the present situation, we can apply Corollary \ref{cor4.2} for the sequences $\rho_m$, $m\in\N$ (previously considered),
\[
h_m:=(a^{(m)}_{-k_m},\ld,a^{(m)}_{-1},a^{(m)}_0), \ \ m\in\N, \ \ \nu_m:=\f_{(a^{(m)}_+,t_m)}\ti, \ \ m\in\N
\]
and $K$, (where $a^{(m)}_+=(a_1^{(m)},\ld,a^{(m)}_{k_m})$ and $a^{(m)}_-=(a^{(m)}_{-k_m},\ld,a^{(m)}_0)(=h_m)$) and find an $m_{\el_0}$ such that, setting $m_{\el_0}=m_0$ for brevity in the notation
\[
\int_X\chi_K\circ\Fi_{a_-^{(m_0)}}d\rho_{m_0}(\nu_{m_0})>1-\oo.
\]
Since $\rho_{m_0}$ is a convex combination, there exists a $z^\ast_{m_0}\in F_{m_0}$ such that
\[
\int_X\chi_K\circ\Fi_{a_-^{(m_0)}}d\f_{z^\ast_{m_0}}(\nu_{m_0})>1-\oo,
\]
i.e. by the form of $\nu_{m_0}$
\begin{eqnarray}
\int_X\chi_K\circ\Fi_{(a_+^{(m_0)},t_{m_0},z^\ast_{m_0},a_-^{(m_0)})}d\ti>1-\oo. \label{eq11}
\end{eqnarray}

We set
\[
\xi_k:=\Fi_{(a_+^{(m_0)},t_{m_0},z^\ast_{m_0},a_-^{(m_0)})}(x_k), \quad k=1,\ld,s
\]
and since
\[
\f_{(a_+^{(m_0)},t_{m_0},z^\ast_{m_0},a_-^{(m_0)})}\bigg(\sum^s_{k=1}\ga_k\de_{x_k}\bigg)=
\sum^s_{k=1}\ga_k\de_{\xi_k}
\]
setting $\cp:=\{k\in\{1,2,\ld,s\}|\xi_k\in K\}$, by (\ref{eq11}) we have
\[
\sum_{k\in\cp}\ga_k>1-\oo.
\]
So, by the definition of the $\ga_k$'s
\[
\sum_{k\in\cp}\nu_{y_0}(Q_k)>(1-\oo)\cdot\nu_{y_0}\bigg(\bigcup^s_{i=1}Q_i\bigg)
\]
and since by (iii) of Lemma \ref{lem5.3} $\nu_{y_0}\Big(\dis\bigcup^s_{i=1}Q_i\Big)>1-\thi$ we have
\begin{eqnarray}
\sum_{k\in\cp}\nu_{y_0}(Q_k)>(1-\oo)(1-\thi).  \label{eq12}
\end{eqnarray}
\noindent
{\bf Claim.} $\Big(\Fi_{(a_+^{(m_0)},t_{m_0},z^\ast_{m_0},a_-^{(m_0)})}\Big(\dis\bigcup_{k\in\cp}\overline{Q}_k\Big)\Big)\cap
\Big(\dis\bigcup^s_{k=1}Q_k\Big)=\emptyset$.

Indeed, by (ii) of Lemma 5.3, diameter${(\overline{Q_k}})=$\,diameter $(Q_k)<h(e)$, for
$k=1,2,\ld,s$, where
$e:=\,$distance$\Big(K,\bcup^s_{k=1}Q_k\Big)$, so we have
\[
\text{diameter}(\Fi_{(a_+^{(m_0)},t_{m_0},z^\ast_{m_0},a_-^{(m_0)})}(\overline{Q_k}))<e, \ \
\text{for} \ \ k=1,2,\ld,s.
\]
On the other hand by the definition of $\cp$, we have
$\xi_k:=\Fi_{(a_+^{(m_0)},t_{m_0},z^\ast_{m_0},a_-^{(m_0)})}(x_k)\in K$, for $k\in\cp$,
where $x_k\in Q_k$. So for $k\in\cp$
\[
(\Fi_{(a_+^{(m_0)},t_{m_0},z^\ast_{m_0},a_-^{(m_0)})}(\overline{Q_k}))\cap\bigg(\bigcup^s_{k=1}Q_k\bigg)=\emptyset
\]
i.e. the claim.

Next, we set
\[
W^\ast:=[\widetilde{pr}^{-1}\{{(a^{(m_0)},t_{m_0},z^\ast_{m_0},a^{(m_0)})}\cap T^{n_{m_0}+|z^\ast_{m_0}|}(B)]\cap[\widetilde{pr}^{-1}\{a\}\cap B].
\]
(where $|z^\ast_{m_0}|$ denotes the length of $z^\ast_{m_0}$)

By (\ref{eq10}) we have $\la(W^\ast)>0$. Clearly, $T^{-(n_{m_0}+|z^\ast_{m_0}|)}(W^\ast)\subseteq B$, so by (\ref{eq12}) and Corollary \ref{cor5.1} we have
\begin{align}
\nu\bigg(\bigg(\bigcup_{k\in\cp}\overline{Q}_k\bigg)\!\!\times\! T^{-(n_{m_0}+|z^\ast_{m_0}|)}(W^\ast)\bigg)&\!>\!((1-\oo)(1-\thi)-\thi)\!\cdot\!\la(T^{-(n_{m_0}+|z^\ast_{m_0}|)}(W^\ast))\nonumber\\
&\!=\!((1\!-\oo)(1\!-\!\thi)-\thi)\cdot\la(W^\ast).  \label{eq13}
\end{align}
Clearly, by the form of $W^\ast$ we have
\begin{equation}
\vPsi^{n_{m_0}+|z^\ast_{m_0}|}\bigg(\!\!\bigg(\bigcup_{k\in\cp}\overline{Q}_k\bigg)\!\times
T^{-(n_{m_0}+|z^\ast_{m_0}|)}(W^\ast)\!\!\bigg)\!=\!(\Fi_{(a_+^{(m_0)},t_{m_0},z^\ast_{m_0},a_-^{(m_0)})}
\bigg(\!\!\bigcup_{k\in\cp}\overline{Q}_k\bigg)\!\!\bigg)\!\times W^\ast  \label{eq14}
\end{equation}
which is measurable, since $\overline{Q}_k$ are compact sets.

By the invariance of $\nu$ under $\vPsi$ and (\ref{eq13}) we have
\begin{align}
\nu\bigg[\vPsi^{n_{m_0}+|z^\ast_{m_0}|}\bigg(\!\!\bigg(\bigcup_{k\in\cp}\overline{Q}_k\bigg)\!\!\times\!
T^{-(n_{m_0}+|z^\ast_{m_0}|)}(W^\ast)\!\bigg)\!\bigg]&\!>\!\nu\!\bigg[\!\!\bigg(\bigcup_{k\in\cp}
\overline{Q}_k\bigg)\!\!\times\!T^{-(n_{m_0}+|z^\ast_{m_0}|)}(W^\ast)\!\bigg] \nonumber \\
&\!>\!((1-\oo)(1-\thi)-\thi)\!\cdot\!\la(W^\ast).  \label{eq15}
\end{align}
By (\ref{eq14}) and (\ref{eq15}) we have
\begin{eqnarray}
\nu\bigg[\bigg(\Fi_{(a_+^{(m_0)},t_{m_0},z^\ast_{m_0},a_-^{(m_0)})}\bigg(\bigcup_{k\in\cp}\overline{Q}_k\bigg)\!\!\bigg)\!\!
\times\!W^\ast\bigg]\!>\!((1-\oo)(1-\thi)-\thi)\!\cdot\!\la(W^\ast).  \label{eq16}
\end{eqnarray}
On the other hand, since clearly $W^\ast\subseteq B$, by (iii) of
Lemma \ref{lem5.3} we have
$\nu_y\Big(\bcup^s_{k=1}Q_k\Big)>1-\thi$, for every $y\in W^\ast$
and intergrating the above inequality over $W^\ast$, we have
\begin{eqnarray}
\nu\bigg(\bigg(\bigcup^s_{k=1}Q_k\bigg)\times
W^\ast\bigg)>(1-\thi)\cdot\la(W^\ast).  \label{eq17}
\end{eqnarray}
Finally, (\ref{eq16}), (\ref{eq17}) and the claim give
\[
\nu(X\times W^\ast)>\frac{3}{2}\cdot\la(W^\ast)
\]
which obviously contradicts the fact that the projection of $\nu$ on
$Y$ coincides with $\la$. \hfill$\square$ \medskip

Finally, combining Theorems \ref{thm6.1} and \ref{thm6.2}, we can conclude the proof of Theorem \ref{thm2.1}. For more details, see [5, Section 8].\vspace*{0.2cm} \\
\noindent
{\footnotesize{\bf Acknowledgements}}. {\footnotesize{We would like to express our gratitude to Professor Constantinos Gryllakis for his guidance during the preparation of this manuscript.}}
%


\begin{thebibliography}{99}
%
\bibitem{1} Brin, M., Stuck, G.: Introduction to Dynamical Systems. Cambridge University Press, Cambridge, 2002.
%
\bibitem{2} Day, M. M.: Fixed point theorems for compact convex sets. Ill. J. Math. {\bf 5} (1961), 585-590.
%
\bibitem{3} Dugundji, J.: Topology, Allyn and Bacon, Boston, 1972.
   %
\bibitem{4} Georgopoulos, P., Gryllakis, C.: Invariant measures for skew products and uniformly distributed sequences. Monatsh. Math. {\bf 167} (2012), no. 1, 81-103.doi:10.1007/s00605-012-0383-z.
%
\bibitem{5} Georgopoulos, P., Gryllakis, C.: Invariant measures for skew products and uniformly distributed sequences II. Monatsh. Math. {\bf 178} (2015), no. 2, 191-220.doi:10.1007/s00605-015-0807-7.
%
\bibitem{6} Hewitt, E., Ross, K. A.: Abstract Harmonic Analysis I, Springer, New York, 1970.
%
\bibitem{7} Kuipers, L., Niederreiter, H.: Uniform Distribution of Sequences. Dover, Mineola, 2006.
    %
\bibitem{8} Mauldin, R. D. (ed.): The Scottish Book. Birkh\"{a}user, Boston, 1981.
%
\bibitem{9} Parthasarathy, K. R.: Probability Measures on Metric Spaces, Academic Press, New York, 1967.
%
\bibitem{10} Ryll - Nardzewski, C.: On the ergodic theorems III. The random ergodic theorem. Studia Math. {\bf 14} (1954), no. 2, 298-301.
%
\bibitem{11} Veech, W.: Some questions of uniform distribution, Ann. Math. (2) {\bf94} (1971), no. 1, 125-138.



\end{thebibliography}
\end{document}